\documentclass[12pt,a4paper,english]{article}
\usepackage[T1]{fontenc}
\usepackage{amsxtra,amssymb,amsmath,enumerate,graphicx,bbm}
\usepackage{theorem}
\usepackage[english]{babel}
\allowdisplaybreaks

\newtheorem{lemma}{Lemma}[section]
\newtheorem{proposition}[lemma]{Proposition}
\newtheorem{theorem}[lemma]{Theorem}
\newtheorem{corollary}[lemma]{Corollary}

\newtheorem{definition}[lemma]{Definition}

\newcommand{\one}{{\rm 1\hspace{-0.1cm}I}}
\newcommand{\PP}{{\mathbbm{P}}}
\newcommand{\DD}{{\mathbb{D}}}
\newcommand{\R}{{\mathbb{R}}}

\newcommand{\E}{{\mathbb{E}}}
\newcommand{\ud}{{\mathrm{d}}}

\newcommand{\mass}{\mathbb{P}}
\newcommand{\m}{\mathbbm{m}}
\DeclareMathAlphabet{\mathitbf}{OML}{cmm}{b}{i}   
\newcommand{\D}{{\mathitbf{D}}}

\newcommand{\DW}{{\DD_{1,2}^0}}
\newcommand{\DN}{{\DD_{1,2}^J}}

\newcommand{\de}{{\delta}}
\newcommand{\om}{{\omega}}
\newcommand{\Om}{{\Omega}}

\newcommand{\s}{{\sigma}}
\newcommand{\al}{{\alpha}}

\newcommand{\ga}{{\gamma}}

\newcommand{\brac}{\left(}
\newcommand{\kets}{\right)}

\newcommand{\kla}{\left ( }
\newcommand{\mer}{\right ) }

\newcommand{\seq}{\subseteq}

\newcommand{\equa}{\begin{eqnarray}}
\newcommand{\tion}{\end{eqnarray}}
\newcommand{\boxx}{\hfill$\square$}
\newcommand{\proof}{{\parindent0pt{\emph{Proof.}}} }

\newcommand{\bor}{\mathcal{B}}

\newcommand{\ftn}{{\cal F}}

\newcommand{\half}{\frac{1}{2}}
\newcommand{\qua}{\begin{eqnarray*}}
\newcommand{\tio}{\end{eqnarray*}}
\newcommand{\pl}{\ \ }                   
\newcommand{\pll}{\ \ \ \ }              
\newcommand{\les}{\hspace{-1.8em}}       
\newcommand{\less}{\hspace{-3.8em}}      
\newcommand{\vare}{\varepsilon}

\begin{document}

\title{Denseness of certain smooth L\'evy functionals in $\DD_{1,2}$}
\author{Christel Geiss \and Eija Laukkarinen}

\maketitle

\begin{center}
Department of Mathematics and Statistics\\
        University of Jyv\"{a}skyl\"{a} \\
        P.O. Box 35 (MaD) \\
        FIN-40014 Jyv\"{a}skyl\"{a} \\
        Finland \\
        chgeiss@maths.jyu.fi \,\,\,\,\,\,
        eija.laukkarinen@jyu.fi
\end{center}


\begin{abstract}
The Malliavin derivative for a L\'evy process $(X_t)$  can be defined 
on the space $\DD_{1,2}$ using a chaos expansion  
or in the case of a pure jump process also via an 
increment quotient operator \cite{sole-utzet-vives}.
In this paper we define the Malliavin  derivative operator $\D$
on the class  $\mathcal{S}$ of smooth random variables 
$f(X_{t_1}, \ldots, X_{t_n}),$ where $f$ is a smooth function with compact 
support. We show that the  closure of 
$L_2(\Om) \supseteq \mathcal{S} \stackrel{\D}{\to} L_2(\m\otimes \mass)$ 
yields to the space 
$\DD_{1,2}.$    
As an application we conclude that 
Lipschitz functions map from $\DD_{1,2}$ into $\DD_{1,2}.$ 
\end{abstract}


\section{Introduction}
In the recent years Malliavin calculus for L\'evy processes has been developed
using various types of chaos expansions. For example, Lee and Shih 
\cite{leh-shih} applied a white noise approach, Le\'on et al. 
\cite{leon-sole-utzet-vives} worked with certain strongly orthogonal 
martingales, L{\o}kka \cite{lokka} and Di Nunno et al. 
\cite{dinunno-meyerbrandis-oksendal-proske} considered multiple integrals with 
respect to the compensated Poisson random measure and Sol\'e et al. 
\cite{sole-utzet-vives2} used the chaos expansions  proved by It\^o \cite{ito}.
  
This chaos representation from It\^o 
applies to any square integrable functional of a general  L\'evy process.  
It uses multiple integrals like in the well known Brownian motion case but
with respect to an independent random measure associated with a L\'evy process.
Sol\'e et al. propose in \cite{sole-utzet-vives} a canonical space for a general 
L\'evy process. They define for random variables on the canonical space 
the increment quotient operator
\[
         \Psi_{t,x}F(\om) =\frac{F(\om_{t,x})-F(\om)}{x}, 
       \hspace{3em}  x \not = 0,
\]
in a pathwise sense, where, roughly speaking, $\om_{t,x}$ can be interpreted 
as the outcome of adding at time $t$ a jump of the size $x$ to the path $\om.$ 
They show that on the canonical L\'evy space the Malliavin derivative $D_{t,x}F$ 
defined via the chaos expansion due to It\^o and  $ \Psi_{t,x}F$
coincide a.e. on $\R_+ \times \R_0 \times \Om$  (where $\R_0:=\R \setminus{0}$)
whenever 
$F \in L_2$ and $\E \int_{\R_+ \times \R_0} |\Psi_{t,x}F|^2 d\m(t,x)< \infty$ 
(see Section 2).\\

On the other hand, on the Wiener space, the Malliavin derivative  is 
introduced as an operator $D$  mapping smooth random  variables 
of the form $F = f(W(h_1), \ldots, W(h_n))$ into $L_2(\Om;H),$ i.e.
\[
  DF= \sum_{i=1}^n \frac{\partial}{\partial x_i} f(W(h_1), \ldots, W(h_n)) h_i,
\] 
(see, for example,  \cite{nualart}).
Here $f$ is a smooth function mapping from $\R^n $ into $\R$ such that 
all its derivatives have at most polynomial growth, and 
$\{W(h), h \in H\}$ is an isonormal Gaussian family associated with 
a Hilbert space $H.$ The closure of the domain of the operator $D$
is the space $\DD_{1,2}.$

In the present paper we proceed in a similar way for a L\'evy process
$(X_t)_{t \ge 0}.$  We will define a Malliavin derivative operator on
a class of smooth random variables and determine its closure. 
The class of smooth random variables we consider consists of elements
of the form $F= f(X_{t_1}, \ldots, X_{t_n})$ where $f:\R^n \to \R$ 
is a smooth function with compact support.

Thanks to the results of Sol\'e et al. \cite{sole-utzet-vives} about
the canonical L\'evy space one can express  the Malliavin derivative 
$DF \in L_2(\m \otimes \mass),$ defined via chaos expansion, on  certain smooth 
random variables of the form $F= f(X_{t_1}, \ldots, X_{t_n})$ explicitly as a 
two-parameter operator $D_{t,x}$
\qua
     D_{t,x}f( X_{t_1},\ldots, X_{t_n}) 
&=& \sum_{i=1}^n \frac{\partial f}{\partial x_i}(X_{t_1},\ldots,X_{t_n})
                    \one_{[0,t_i]\times \{ 0 \} }(t,x) \\
   &&+ \Psi_{t,x} f( X_{t_1},\ldots, X_{t_n}) \one_{\{x\not =0\}}(x),
\tio
for $\m\otimes\mass$-a.e. $(t,x,\om).$ Here $\Psi_{t,x},$ for $x\not =0,$
is given by
\qua
&& \les \Psi_{t,x} f( X_{t_1},\ldots, X_{t_n})\\ 
&:=& 
\frac{f(X_{t_1} +
x\one_{[0,t_1]}(t),\ldots,X_{t_n} +
x\one_{[0,t_n]}(t))-f(X_{t_1},\ldots,X_{t_n})}{x}. 
\tio
Our main result is that the smooth random variables 
$f( X_{t_1},\ldots, X_{t_n})$ are dense in the space $\DD_{1,2}$ 
defined via the chaos expansion. This implies that
defining $\D$ as an operator on the  smooth random 
variables as in Definition \ref{operator-on-S} below and taking the closure 
leads to the same result as defining  
$D$ using  It\^o's chaos expansion (see Definition \ref{operator-on-chaos}).

The paper is organized as follows. In Section 2 we shortly recall 
It\^o's chaos expansion, the definition of the  Malliavin derivative and
some related facts.
In Section 3 we introduce the Malliavin derivative operator on smooth 
random variables. In Section 4 we determine its closure. Applying the denseness
result from the previous  section we show in Section 5 that Lipschitz functions 
map from  $\DD_{1,2}$ into $\DD_{1,2}.$ \\

\section{The Malliavin derivative  via It\^{o}'s chaos expansion}

We assume  a c\'adl\'ag  L\'evy process $X = (X_t)_{t\ge 0}$ on
a complete probability space $(\Om,\ftn,\PP)$ with L\'evy triplet
$(\ga,\s^2,\nu)$ where $\ga\in\R$, $\s \geq 0$ and $\nu$ is the L\'evy
measure. Let $N$ be the  Poisson random measure associated 
with the process $X$ and $\tilde{N}$ the compensated Poisson random
measure, $d \tilde{N}(t,x) = d N(t,x) - d t d \nu(x)$.

Consider the measures $\mu$ on $\bor(\R)$, 
\qua
d \mu(x) := \s^2 d \de_0(x) + x^2 d \nu(x),
\tio
and $\m$ on $\bor(\R_+\times\R)$, where $\R_+:=[0,\infty),$
\qua
d\m(t,x) := d t d \mu(x).
\tio 
For $B\in \bor(\R_+\times\R)$ such
that $\m(B) < \infty$ let
\qua
M(B) = \s \int_{\{t\in\R_+:(t,0)\in B\}} \ud W_t
+ \lim_{n\to\infty}\int_{\{(t,x)\in B: 1/n < |x| < n\}} x d \tilde{N}(t,x),
\tio
where the convergence is taken in $L_2(\Om,\ftn,\PP)$. Now $\E
M(B_1) M(B_2) = \m(B_1 \cap B_2)$ for all $B_1,B_2$ with
$\m(B_1)<\infty$ and $\m(B_2)<\infty$.
For $n=1,2,\ldots$ write 
\qua
L^n_2 := L_2 \brac (\R_+\times\R)^n,
\bor(\R_+\times\R)^{\otimes n}, \m^{\otimes n} \kets.
\tio
For $f\in L^n_2$ It\^{o} \cite{ito} defines a
multiple integral $I_n(f)$ with respect to the random measure $M$.
It holds $I_n(f) = I_n(\tilde{f}),$ a.s., where $\tilde{f}$ is the
symmetrization of $f$,
$$\tilde{f}(z_1,\ldots,z_n) 
= \frac{1}{n!} \sum_{\pi_n} f(z_{\pi(1)},\ldots,z_{\pi(n)})
\quad \textrm{for all } z_i=(t_i,x_i) \in\R_+\times\R.$$
The sum is taken over all permutations
$\pi_n:\{1,\ldots,n\}\to\{1,\ldots,n\}$. 

Let $(\ftn^X_t)_{t\geq0}$ be  the augmented natural filtration of $X$.
Then $(\ftn_t^X)_{t\geq0}$ is right continuous 
(\cite[Theorem  I 4.31]{protter}). 
Set $\ftn^X := \bigvee_{t\geq 0} \ftn^X_t.$ By Theorem 2 of It\^{o}
\cite{ito} it holds the chaos decomposition 
$$L_2 := L_2(\Om,\ftn^X,\PP) = \bigoplus_{n=0}^\infty I_n(L^n_2), $$
where $I_0(L^0_2):= \R$ and
$I_n(L^n_2):=\{I_n(f_n):f_n\in L^n_2\}$ for $n=1,2,\ldots.$  For
$F \in L_2$ the representation
$$F = \sum_{n=0}^\infty I_n(f_n),$$
with $I_0(f_0) = \E F,$ a.s., is unique if the functions  $f_n$ are symmetric.
Furthermore,
\qua
\|F\|^2_{L_2} = \sum_{n=0}^\infty n! \|\tilde{f}_n\|^2_{L^n_2}. 
\tio
\begin{definition} \label{operator-on-chaos}
Let  $\DD_{1,2}$ be the space of all $F = \sum_{n=0}^\infty
I_n(f_n) \in L_2$ such that
$$\|F\|^2_{\DD_{1,2}} := 
\sum_{n=0}^\infty (n+1)! \|\tilde{f}_n\|^2_{L^n_2} < \infty.$$
Set 
$
L_2(\m\otimes\mass) := L_2(\R_+\times\R\times\Om,
\bor(\R_+\times\R)\otimes\ftn^X, \m\otimes\PP).$ 
The Malliavin derivative operator 
$D:\DD_{1,2}\to L_2(\m\otimes\PP)$  is defined by
\qua
D_{t,x}F := \sum_{n=1}^\infty n I_{n-1} (\tilde{f}_n((t,x),\cdot)), 
 \pll (t,x,\om)\in \R_+\times\R\times\Om.
\tio
\end{definition}
\bigskip
We consider (as Sol\'e et al. \cite{sole-utzet-vives}) the operators 
$D_{\cdot,0}$ and $D_{\cdot,x}, \, x\neq0$ and their domains $\DW$ and 
$\DN$.
For $\s>0$ let $\DW$ be consisting of
random variables $ F = \sum_{n=0}^\infty I_n(f_n) \in L_2$ such
that
$$\|F\|^2_{\DW} := \|F\|^2_{L_2} + \sum_{n=1}^\infty n\cdot n!
\|\tilde{f}_n  \one_{(\R_+\times\{0\})\times(\R_+\times\R)^{n-1}}  \|^2_{L_2^n} <\infty.$$
For $\nu \neq 0$, let $\DN$ be the set of $F\in L_2$ such that
\qua
  \|F\|^2_{\DN} 
:=\|F\|^2_{L_2} +  \sum_{n=1}^{\infty} n\cdot n!
    \|\tilde{f}_n  \one_{(\R_+\times\R_0) \times(\R_+\times\R)^{n-1}}\|^2_
    {L_2^n} 
< \infty,
\tio
where $\R_0:=\R \setminus \{0\}.$
If both $\s > 0$ and $\nu \neq 0$, then it holds
\equa \label{intersection} 
\DD_{1,2} = \DW \cap \DN.
\tion

In case $\nu = 0$,  $D_{\cdot,0}$
coincides with the classical Malliavin derivative $D^W$ (see, for
example,  \cite{nualart}) except for a multiplicative
constant, $D^W_t F = \s D_{t,0} F$.\\ 

In the next  lemma we formulate a denseness result which will be used
to determine the closure of the Malliavin operator from Definition 
\ref{def:malliavin} below. 

\begin{lemma} \label{lemma:span_dense}
Let $\mathcal{L}\seq L_2$ be the linear span of random variables of the form
\qua
     M(T_1\times A_1) \cdots M(T_N \times A_N), \pll N=1,2,\ldots
\tio
where the  $A_i's$ are finite intervals of the form $(a_i,b_i]$  and the
$T_i's$ are finite disjoint 
intervals of the form $T_i =(s_i,t_i]$.
Then $\mathcal{L}$ is dense in $L_2,$ $\DD_{1,2},$ $\DW$
and $\DN.$ 
\end{lemma}
\bigskip
\proof $1^\circ$ First we consider the class of all linear 
combinations of 
\qua 
  M(B_1) \cdots M(B_N)= I_N(\one_{B_1 \times \cdots \times B_n}),
\tio
$N=1,2, \ldots,$ where the sets
$B_i \in \bor(\R_+\times\R)$ are disjoint and fulfill $\m(B_i) < \infty.$
It follows from the completeness of the multiple integrals 
in $L_2$ (see \cite[Theorem 2]{ito}) that this class 
is dense in $L_2.$ Especially, the class of all linear 
combinations of $\one_{B_1 \times \cdots \times B_n}$ with disjoint sets $B_i$ of 
finite measure $\m$ is dense in 
$L^n_2 = L_2 \brac (\R_+\times\R)^n,
\bor(\R_+\times\R)^{\otimes n}, \m^{\otimes n} \kets.$ Let 
$\mathcal{L}_n$ be the linear span of $\one_{T_1\times A_1} \otimes  
\cdots \otimes  \one_{T_n \times A_n}$ with disjoint intervals $T_i.$ 
One can easily see that  $\mathcal{L}_n$ is dense in $L^n_2$ as well:\\
Because $\m$ is a Radon measure, there are compact sets 
$C_i \seq B_i$ such that $\m(B_i \setminus C_i)$ is sufficiently small
to get 
\qua
\| \one_{B_1 \times \cdots \times B_n}
    - \one_{C_1 \times \cdots \times C_n}\|_{L^n_2} <\vare
\tio
for some given $\vare >0.$ Since the compact sets $(C_i)$ are disjoint 
one can find disjoint bounded open sets 
$U_i \supseteq C_i$ such that $\| \one_{C_1 \times \cdots \times C_n}
    - \one_{U_1 \times \cdots \times U_n}\|_{L^n_2} <\vare.$
For any bounded open set $U_i \seq (0, \infty) \times \R$ one can find a 
sequence of 'half-open rectangles' 
$Q_{i,k}= (s_k^i,t_k^i]\times(a_k^i,b_k^i]=T_k^i\times A_k^i,$ such that  
$U_i = \bigcup_{k=1}^{\infty} Q_{i,k}$
(taking  half-open rectangles $Q_x \seq U_i$  with rational 'end points'
containing the point $x \in U_i$ gives
$U_i = \bigcup_{Q_x \seq U_i}^{\infty} Q_x$).

Hence for sufficiently large $K_i's$ one has
\qua
\| \one_{U_1 \times \cdots \times U_n}
    - \one_{ \bigcup_{k=1}^{K_1}
Q_{1,k} \times \cdots \times \bigcup_{k=1}^{K_N}
Q_{n,k}}\|_{L^n_2} <\vare,
\tio
where the $Q_{i,k}$ can now be chosen such that they are disjoint.
This implies that the linear span of 
$\one_{Q_1 \times \cdots \times Q_n},$  where the $Q_i's$
are of the form $T_i\times A_i$ is dense in $L^n_2 .$

$2^\circ$ 
For the convenience of the reader we recall the idea of the proof of 
Lemma 2 \cite{ito} to show that the intervals $T_i$ can be chosen disjoint.
Consider \equa
     M(T_1\times A_1) \cdots M(T_n \times A_n),
\tion 
with $\mu(A_1) >0, \ldots,\mu(A_n) >0,$ 
where $(T_i\times A_i)\cap (T_j\times A_j) = \emptyset$ for
$i\neq j$. Assume, for example (all other cases can be treated similarly), 
that $T_1 = \dots = T_m =:T$ while  $T_m, \ldots, T_n$ are pairwise disjoint.  
Given the expression
\qua
M(T\times A_1)\cdots M(T\times A_m) M(T_{m+1}\times A_{m+1})
 \cdots M(T_n\times A_n)
\tio
choose an equidistant partition $(E_j)_{j=1}^N$ of $T$ so that 
$|E_j| = \frac{|T|}{N}$
and set 
$$c:=(n+1)\mu(A_1)\cdots\mu(A_n) |T_{m+1}|\cdots |T_n|.$$
 Now
\qua
&& \les   M(T\times A_1)\cdots M(T\times A_m) M(T_{m+1}\times
    A_{m+1})\cdots M(T_n\times A_n)\\
& = &\sum_{\stackrel{j_1, \ldots, j_m =1}
     {\text{ all }j_i \text{ distinct}}}^N\hspace{-1em} M(E_{j_1}\times
    A_1)\cdots  M(E_{j_m}\times A_m)\\
&& \hspace{8em} \times    M(T_{m+1}\times
    A_{m+1})\cdots M(T_n\times A_n)\\
&& + \sum_{\stackrel{j_1, \ldots, j_m =1}
     {j_i \text{ not distinct}}}^N
   \hspace{-1em} M(E_{j_1}\times A_1)\cdots
    M(E_{j_m}\times A_m) \\
&& \hspace{8em} \times  M(T_{m+1}\times A_{m+1})\cdots M(T_n\times A_n)\\
& =& S_1 + S_2,
\tio
where $S_1$ is a sum of products with
disjoint time intervals. We complete the proof by observing that
\qua
    \|S_2\|^2_{\DD_{1,2}}
&=& (n+1) \hspace{-1em}\sum_{\stackrel{j_1, \ldots, j_m =1}
     {j_i \text{ not distinct}}}^N \hspace{-1em}
    \m(E_{j_1}\times A_1)\cdots
    \m(E_{j_m}\times A_m) \\
&&  \hspace{10em} \times \m(T_{m+1}\times A_{m+1})
    \cdots \m(T_n\times A_n)\\
&=& c \sum_{\stackrel{j_1, \ldots, j_m =1}
     {j_i \text{ not distinct}}}^N
    |E_{j_1}|\cdots|E_{j_m}|\\
& =& c \kla N^m-N (N-1) \cdots (N-(m-1))\mer \kla\frac{|T|}{N}\mer^m \\
&=&  c|T|^m (1-(1-\frac{1}{N})\cdots (1-\frac{m-1}{N}))  \to 0
\tio
for $N \to \infty.$\\

$3^\circ$ 
The denseness of  $\mathcal{L}_n$ in $L^n_2$ implies that    $\mathcal{L}$
is dense in $L_2$ and $\DD_{1,2}.$ The remaining cases follow from the 
fact that 
\qua
\|f_n \one_{(\R_+\times\{0\})\times(\R_+\times\R)^{n-1} } \|_{L_2^n} 
\le \|f_n\|_{L_2^n}
\tio
and
\qua
   \|f_n\one_{(\R_+\times\R_0)\times(\R_+\times\R)^{n-1},\, \m^{\otimes n}) }  
   \|_{L_2^n}
\le \|f_n\|_{L_2^n}.
\tio
\boxx


\section{The Malliavin derivative as operator on  $\mathcal{S}$}
Let $C_c^\infty(\R^n)$ denote the space of smooth functions
$f: \R^n \to \R$ with compact support. 
\begin{definition} \label{def:malliavin}
A random variable of the form  $F=f(X_{t_1},\ldots,X_{t_n}),$ where 
$f\in C_c^\infty(\R^n)$, $n\in\mathbb{N},$ and  $t_1,\ldots,t_n\geq0$
is said to be a smooth random variable. The set of all
smooth random variables is denoted by $\mathcal{S}.$ 
\end{definition}

\begin{definition} \label{operator-on-S}
For  $F=f(X_{t_1},\ldots,X_{t_n}) \in \mathcal{S} $ we
define the Malliavin derivative operator $\D$  as a map from 
$\mathcal{S}$ into $L_2(\m\otimes\PP)$ by
\qua
 && \les \D_{t,x} f(X_{t_1},\ldots,X_{t_n})\\
 &:=& 
  \sum_{i=1}^n \frac{\partial f}{\partial x_i}
          (X_{t_1},\ldots,X_{t_n})\one_{[0,t_i]\times\{0\}}(t,x)\\
\les &&\les + \frac{f(X_{t_1} \! + x\one_{[0,t_1]}(t),\ldots,X_{t_n} \! +
x\one_{[0,t_n]}(t))-f(X_{t_1},\ldots,X_{t_n})}{x}\one_{\R_0}(x)
\tio 
for  $(t,x) \in \R_+ \times \R.$
\end{definition}
\bigskip
Since for $f(X_{t_1},\ldots, X_{t_n}) \in \mathcal{S}$ it holds
\qua
\E \int_{\R_+}|\D_{t,0} f(X_{t_1},\ldots,X_{t_n})|^2 dt < \infty,
\tio
one can apply the methods used in the proof of Proposition 3.5 
in \cite{sole-utzet-vives} for the canonical L\'evy space to show that 
\equa  \label{brownian-connection}
\D_{t,0} f(X_{t_1},\ldots,X_{t_n}) 
&=&  D_{t,0} f(X_{t_1},\ldots,X_{t_n}), \pl dt\otimes \mass -\text{a.e.}
\tion
We conclude $\mathcal{S} \seq \DW.$ Because we have also  
\qua
 \E \int_{\R_+\times\R_0} |\D_{t,x} f(X_{t_1},\ldots,X_{t_n})|^2 d\m(t,x) < \infty
\tio
for  all $f(X_{t_1},\ldots,X_{t_n}) \in \mathcal{S}$
one can show by the same methods as for \eqref{brownian-connection} that
\qua \D_{t,x} f(X_{t_1},\ldots,X_{t_n}) 
=   D_{t,x} f(X_{t_1},\ldots,X_{t_n})  
\tio
for $\m \otimes \mass$ a.e. $(t,x,\om)
 \in \R_+\times\R_0 \times \Om$ and therefore $\mathcal{S} \seq \DN.$ 
See \cite[Proposition 5.5]{sole-utzet-vives} for an according relation 
on the canonical L\'evy space. Consequently:

\begin{lemma} It holds
$\D F=DF$ in $L_2(\m \otimes \mass)$  for all
$F \in \mathcal{S}.$
\end{lemma}
 Especially,
we conclude from the fact that any $F \in L_2 \supseteq \mathcal{S}$ has a 
unique chaos expansion 
that also  $\D F$ does not depend from the  representation 
$F=f(X_{t_1},\ldots,X_{t_n}) \in \mathcal{S}.$ \\

Using the equality of $\D$ and $D$ on $\mathcal{S}$ and the fact
that $\mathcal{S}$ is closed with respect to multiplication
we are now able to reformulate Proposition 5.1 of \cite{sole-utzet-vives}
for our situation:  

\begin{corollary} 
\label{corollary:not-productrule}
For $F$ and $G$  in $\mathcal{S}$ it holds 
\qua
D_{t,x}(FG) = GD_{t,x}F+FD_{t,x}G +x D_{t,x}FD_{t,x}G
\tio
for $\m\otimes\mass-$a.e. $(t,x,\om) \in \R_+ \times \R \times \Om.$
\end{corollary}
\bigskip


\section{The closure of the Malliavin derivative operator}

The operator $\D: \mathcal{S} \to  L_2(\m \otimes \mass)$ is closable, if for 
any sequence $(F_n) \seq \mathcal{S}$ which converges to $0$
in $L_2$  such that  $\D(F_n)$  converges in $ L_2(\m \otimes \mass),$ 
it follows that  $(\D F_n)$ converges to $0$ in  $L_2(\m \otimes \mass).$ 
As we know from the previous section that $D$ and $\D$ coincide on $\mathcal{S} 
\seq \DD_{1,2},$ it is clear that $\D$ is closable and the closure of the domain of 
definition of  $\D$ with respect to the norm 
\qua
\|F\|_{\D} := \big [ \E |F|^2 + \E \|\D F \|^2_{L_2(\m)} \big ]^{\half},
\tio
is  contained in  $\DD_{1,2}.$ What remains to show is that the closure is 
equal to  $\DD_{1,2}.$

\begin{theorem}\label{theorem:closure=DD_{1,2}}
The closure of $\mathcal{S}$ with respect to the norm 
$\|\cdot\|_{\D} $
is the space $\DD_{1,2}.$
\end{theorem}
\bigskip
Theorem \ref{theorem:closure=DD_{1,2}} implies that 
the  Malliavin derivative $D$ defined via It\^o's chaos expansion and the 
closure of the operator 
$L_2 \supseteq \mathcal{S} \stackrel{\D}{\to} L_2(\m\otimes \mass)$ 
coincide. Before we start with the proof we formulate a Lemma
for later use.
%
%
\begin{lemma} \label{lemma:first-chaos-approximation}
For $\varphi \in C_c^{\infty}(\R)$ and partitions 
$\pi_n:=\{s=t^n_0 <t^n_1<\dots<t^n_n=u\}$ of the interval $[s,u]$  it 
holds for $\psi(x) := x\varphi(x)$ that
\qua
\DD_{1,2}-\lim_{|\pi_n|\to 0} \bigg (\sum_{j=1}^n \psi(X_{t^n_j} - X_{t^n_{j-1}}) -
\E&&\les \sum_{j=1}^n \psi(X_{t^n_j} - X_{t^n_{j-1}}) \bigg ) \\
&=& \int_{(s,u]\times\R} \varphi(x)\ d M(t,x),
\tio
where $|\pi_n|:=\max_{1\le i \le n} |t_i^n-t_{i-1}^n|.$
\end{lemma}
\proof
To keep the notation simple, we will drop the $n$ of the partition 
points $t^n_j.$ Notice that 
$\int_{(s,u]\times\R} \varphi(x)\ d M(t,x)= I_1(\one_{(s,u]}\otimes\varphi).$
We set
\qua
G^n:=  \sum_{j=1}^n \psi(X_{t_j} - X_{t_{j-1}}) -
\E\sum_{j=1}^n \psi(X_{t_j} - X_{t_{j-1}})
\tio
and
\qua
G:=\int_{(s,u]\times\R} \varphi(x)\ d M(t,x).
\tio
We write $D_{t,x} G^n$  explicitly as
\qua
&& \less  D_{t,x} \kla \sum_{j=1}^n \psi(X_{t_j} - X_{t_{j-1}}) \mer \\
&=& \sum_{j=1}^n \psi'(X_{t_j} - X_{t_{j-1}}) 
       \one_{(t_{j-1},t_j] \times\{0\}}(t,x) \\
&&+ \sum_{j=1}^n \frac{\psi(X_{t_j} - X_{t_{j-1}}+x)-
       \psi(X_{t_j} - X_{t_{j-1}})}{x} \one_{(t_{j-1},t_j]\times\R_0}(t,x).
\tio
Moreover, we have 
$D_{t,x}I_1(\one_{(s,u]}\otimes\varphi)=\one_{(s,u]}(t)\varphi(x)$ $\m-$a.e. 
Using the general fact that for any $F\in L_2$ with expectation zero it 
holds  $\|F\|^2_{\DD_{1,2}} \le 2 \|DF\|^2_{L_2(m\otimes\PP)}$
we obtain  
\qua
&&  \les  \|G-G^n\|^2_{\DD_{1,2}} \\
& \le & 2 \|DG- DG^n\|^2_{L_2(\m\otimes\PP)}\\
& =& 2 \E \int_{\R_+} \sum_{j=1}^n \one_{(t_{j-1},t_j]}(t)
    \big  [\varphi(0)-\psi'(X_{t_j} - X_{t_{j-1}}) \big ]^2 dt\\
&& +2\E \int_{\R_+\times \R_0}\sum_{j=1}^n \one_{(t_{j-1},t_j]}(t) 
 \big  [\varphi(x)- D_{t,x}\psi(X_{t_j} - X_{t_{j-1}})\big ]^2  d \m(t,x)\\
& =&  2\E \int_{\R_+} \sum_{j=1}^n \one_{(t_{j-1},t_j]}(t)
    \big  [\varphi(0)-\psi'(X_{t_j} - X_{t_{j-1}}) \big ]^2 dt
        \\
&&+ 2\E \int_{\R_+\times \R_0}\sum_{j=1}^n   \one_{(t_{j-1},t_j]} (t)
     [ \psi(X_{t_j}- X_{t_{j-1}} + x) \\
&&    \hspace{8em}    - \psi(X_{t_j} - X_{t_{j-1}})- \psi(x) ]^2 dtd\nu(x)\\
&& \to 0
\tio
as $n\to\infty$ because of dominated convergence and the a.s. c\`adl\`ag
property of the paths of $(X_t).$ Indeed, one can use the estimates
\qua
|\varphi(0) - \psi'(y)|\le\|\varphi\|_\infty+ \|\psi'\|_\infty,
\tio
\qua
|\psi(y + x) -
\psi(y) - \psi(x)| \leq (\|\psi'\|_\infty+\|\varphi\|_\infty  
  + 3\|\psi\|_\infty)(|x| \wedge 1),
\tio 
and that $\int_{\R} (x^2\wedge 1) d\nu(x) < \infty.$
Moreover, for $|\pi_n| \to 0$ we have  from the  c\`adl\`ag
property of the paths the pointwise convergence in $t \in (s,u]$ of
\qua
\sum_{j=1}^n \one_{(t^n_{j-1},t^n_j]}(t)
    \big  [\varphi(0)-\psi'(X_{t^n_j} - X_{t^n_{j-1}}) \big ]^2 \to
[\varphi(0)-\psi'(X_t - X_{t-}) \big ]^2
\tio 
(note that $\varphi(0)- \psi'(0)=0$), and of
\qua
&& \less \sum_{j=1}^n   \one_{(t^n_{j-1},t^n_j]} (t)
     [ \psi(X_{t^n_j}- X_{t^n_{j-1}} + x) 
     - \psi(X_t - X_{t-})- \psi(x) ]^2 \\
&& \to [ \psi(X_t- X_{t-} + x) 
     - \psi(X_t - X_{t-})- \psi(x) ]^2.
\tio
Because the set $\{ t>0; \pl X_t - X_{t-} \not=0\}$ 
is at most countable  for c{\`a}dl{\`a}g paths the assertion follows. 
\boxx \\
\bigskip \\
{\it Proof of Theorem \ref{theorem:closure=DD_{1,2}}}.
According to Lemma \ref{lemma:span_dense} it is sufficient to show that
 an expression like
$ M(T_1\times A_1) \cdots M(T_N \times A_N),$
where the $A_i's$ are bounded Borel sets  and the $T_i's$ finite  disjoint 
intervals, can be approximated in  $\DD_{1,2}$ by a sequence 
$(F_n)_n \seq \mathcal{S}.$  \\

$1^{\circ}$  In this  step we want to show that it is enough to 
approximate 
\equa \label{smooth-product}
I_1(\one_{T_1}\otimes\varphi_1)\cdots I_1(\one_{T_N} \otimes\varphi_N),
\tion
by $(F_n)_n \seq \mathcal{S}$ where $\varphi_{i} \in  C_c^\infty(\R).$ 
Since the intervals $T_i$ are disjoint the definition of the 
multiple integral implies that 
\qua
 M(T_1\times A_1) \cdots M(T_N \times A_N)
= I_N(\one_{T_1\times A_1} \otimes \cdots \otimes 
  \one_{T_N \times A_N}) \pll \text{a.s.}
\tio
By the same reason,
\qua
I_1(\one_{T_1}\otimes \varphi_1) \cdots  I_1(\one_{T_N} \otimes \varphi_N)
= I_N((\one_{T_1}\otimes \varphi_1) \otimes \cdots \otimes
 (\one_{T_N}\otimes \varphi_N))   \pll \text{a.s.}   
\tio
We have 
\qua
&& \| I_N(\one_{(T_1\times A_1)\times \cdots \times (T_N \times A_N)})-
I_N((\one_{T_1}\otimes \varphi_1) \otimes \cdots \otimes
 (\one_{T_N}\otimes \varphi_N) )\|^2_{\DD_{1,2}} \\
&\le& (N+1)!\|\one_{(T_1\times A_1)\times \cdots \times (T_N \times A_N)} -  
    (\one_{T_1}\otimes \varphi_1) \otimes \cdots \otimes
 (\one_{T_N}\otimes \varphi_N)       \|^2_{L_2^N} \\
&\le& (N+1)! |T_1| \cdots |T_N| 
 \|\one_{ A_1 \times \cdots \times A_N}- \varphi_1 \otimes\cdots\otimes 
 \varphi_N\|^2_{L_2^N(\mu^{\otimes N})}.
\tio
The last expression can be made arbitrarily  small by choosing  $\varphi_i$
such that $\|\one_{ A_i}-\varphi_i\|_{L^1_2(\mu)}$ is small. For example,
for each $i$ there are compact sets 
$C^i_1 \subseteq C^i_2 \subseteq \cdots
\subseteq A_i$  and open sets $U^i_1 \supseteq U^i_2 \supseteq
\cdots \supseteq A_i$ such that 
\qua
           \mu(U^i_n \setminus C^i_n) \to 0
\tio 
as $n\to\infty$. By the $C^\infty$ Urysohn
Lemma (\cite{folland}, p. 237) there is for each $n$ a
function $\varphi^i_n\in C_c^\infty(\R)$ such that
$0\le \varphi^i_n \le 1$, $\varphi^i_n=1$ on $C^i_n$ and
supp$(\varphi^i_n) \subset U^i_n$. Then
\qua
\|\one_{ A_i}-\varphi^n_i\|^2_{L^1_2(\mu)} \le  \mu(U^i_n \setminus C^i_n)  \to 0 
\tio
as $n \to \infty.$ \\

$2^{\circ}$ Now we use Lemma \ref{lemma:first-chaos-approximation}
to approximate the  expression \eqref{smooth-product} by a sequence  
$(F_n)_n \seq \mathcal{S}.$
For $i=1,\ldots,N$ set $\psi_i(x):=x\varphi_i(x)$ and
\qua
G_i^n :=\sum_{j=1}^n \one_{\{t_j,t_{j-1} \in \bar{T}_i\}}
   \psi_i(X_{t_j} - X_{t_{j-1}}) -
\E\sum_{j=1}^n\one_{\{t_j,t_{j-1} \in \bar{T}_i\}}\psi_i(X_{t_j} 
 - X_{t_{j-1}}).
\tio
The partition $\pi_n$ can be chosen such that 
all end points of the intervals $\bar{T}_i$ belong to  $\pi_n.$
Let
\qua
G_i := I_1(\one_{T_i}\otimes \varphi_i).
\tio 
Because the intervals $(T_i)$ are disjoint it follows 
from Corollary \ref{corollary:not-productrule} that the product rule holds
in our case:
\equa \label{product-rule}
D\Pi_{i=1}^N G_i^n = \sum_{i=1}^N G^n_1 \cdots G_{i-1}^n(DG_i^n) G_{i+1}^n
           \cdots  G_N^n \pll \m\otimes \mass-\text{a.e.}  
\tion
Indeed, because  of $D_{t,x}G_i^n = (D_{t,x}G_i^n)\one_{T_i}(t)$ it follows
\qua
x(D_{t,x}G_{i_1}^n)\one_{T_{i_1}}(t) (D_{t,x}G_{i_2}^n)\one_{T_{i_2}}(t) =0
\pll \m\otimes \mass-\text{a.e.}  
\tio
for any $i_1 \not =i_2$. Equation \eqref{product-rule} follows then by 
induction. 

We observe that $ G^n_1,  \ldots,  G_N^n$ as well as 
$G^n_1, \ldots, G_{i-1}^n, DG_i^n, G_{i+1}^n, \ldots,  G_N^n$
are  mutually independent by construction. Hence to show $L_2-$convergence
of these products it is enough to have  $L_2-$convergence for each factor. From
Lemma \ref{lemma:first-chaos-approximation} we obtain that
$G_i^n \to G_i$ in $\DD_{1,2}$ for all $i=1,\ldots, N,$ so that
\qua
L_2(\m \otimes \mass)-\lim_{|\pi_n| \to 0} 
           G^n_1 \cdots G_{i-1}^n(DG_i^n) G_{i+1}^n
           \cdots  G_N^n \\
=  G_1 \cdots G_{i-1}(DG_i) G_{i+1}
           \cdots  G_N.
\tio
Hence we have found a sequence $f_n(X_{t_0},\ldots,X_{t_n})=\Pi_{i=1}^N G_i^n$
converging to expression \eqref{smooth-product}  
in $\DD_{1,2}$ where $f_n \in
C^\infty(\R^{n+1})$. 

To find a sequence $(F_n) \seq \mathcal{S}$  we  proceed as follows.
Choose  functions $\beta_N \in
C_c^{\infty}(\R)$ such that $0 \le \beta_N\le 1$ and  $\beta_N(x)=1$ for
$|x| \le N,$ the support of  $\beta_N$ is contained in  
$\{x; |x| \le N +2\}$  and 
$\|\beta'_N\|_{\infty}\le1.$ Set $x_{-1}:=0$ and
$\al_N(x_0, \ldots,x_n) := \Pi_{i=0}^n \beta_N(x_i-x_{i-1}).$ 
Then we have $\al_N (x) f_n(x) \in C_c^{\infty}(\R^{n+1}),$ 
and with similar estimates like above one can show that for
$0 \le t_0 \le \ldots \le t_n $ we have
\qua
\al_N (X_{t_0},\ldots,X_{t_n}) f_n(X_{t_0},\ldots,X_{t_n})
\to f_n(X_{t_0},\ldots,X_{t_n})
\tio
in  $\DD_{1,2}$ as $N \to \infty.$   

\boxx\\

\begin{corollary} \label{corollary:denseness}
The set $\mathcal{S}$ of smooth random variables is 
dense in  $\DW$ and $\DN.$
\end{corollary}
\bigskip
\proof Assume $F \in \DW$ has the representation 
$F= \sum_{m=1}^{\infty} I_m(f_m).$  
For a given $\vare >0$ fix $N$ such that 
$\| \sum_{m=N}^{\infty} I_m(f_m)\|_{\DW} < \vare.$ From $F \in L_2$ we 
conclude $F^N := \sum_{m=1}^{N} I_m(f_m) \in \DD_{1,2}.$  
By Theorem \ref{theorem:closure=DD_{1,2}} we can find a sequence 
$(F_n) \seq \mathcal{S}$ converging to $F^N$ in $\DD_{1,2}$
and therefore also in $\DW.$  In the same way one can see that  $\mathcal{S}$ 
is dense in  $\DN.$
\boxx\\

\section{Lipschitz functions operate on  $\DD_{1,2}$}


\begin{lemma} \label{Lipschitz-lemma}
Let $g :\R \to\R$ be  Lipschitz continuous with Lipschitz constant $L_g.$  
\begin{enumerate}[(a)]
\item If $\s>0$, then $g(F) \in   \DW$ for all $F  \in   \DW$ and 
\equa \label{D_{t,0}-part}
     D_{t,0} g(F) = G D_{t,0}F \quad
                    dt \otimes\mass -\text{a.e.},
\tion
where $G$ is a random variable which is a.s. bounded by  $L_g.$
\item If $\nu \not =0,$  then $g(F) \in\DN$ for all $F\in\DN,$ where
\equa \label{D_{t,x}-part} 
   D_{t,x}g(F) 
   = \frac{g(F + xD_{t,x}F)-g(F)}{x}
\tion
for $\m\otimes\PP\text{-a.e.}$ $(t,x,\om) \in \R_+ \times \R_0 \times \Om.$ 
\end{enumerate}
\end{lemma}

\proof We will adapt the proof of Proposition 1.2.4 of \cite{nualart}
to our situation.
Corollary \ref{corollary:denseness} implies that there exists
a sequence $(F_n) \seq  \mathcal{S}$ of the form 
$F_n=f_n( X_{t_1},\ldots,X_{t_n})$ which converges to $F$ in $\DW$.
Like  in \cite{nualart} we choose a non-negative $\psi \in 
C_c^{\infty}(\R)$ such that supp$(\psi) \seq [-1,1]$  and 
$\int_{\R} \psi(x)dx =1$ 
and set  $\psi_N(x):= N\psi(Nx).$ 

Then  $g_N:=g*\psi_N$ is smooth and converges  to $g$ uniformly.
 Moreover, $\|g_N'\|_{\infty} \le L_g.$
Hence  $g_N \circ F_n=g_N(f_n( X_{t_1},\ldots,X_{t_n}))\in \mathcal{S}$
and $(g_n(F_n))$ converges to $g(F)$ in $L_2.$ \\
Moreover, 
\qua
&& \less \E  \int_{\R_+} |D_{t,0} g_n(F_n)|^2 dt \\
&=& \E g_n'(F_n)^2  \int_{\R_+} \big | \sum_{i=1}^n 
      \frac{\partial}{\partial x_i}f_n(X_{t_1}, \ldots, X_{t_n}) 
     \one_{[0,t_i]}(t) \big |^2 dt \\
&\le & L_g^2 \, \|F_n\|_{\DW}^2.
\tio 
Since we have convergence of $(g_n(F_n))$ to $g(F)$ in $L_2$ and
it holds 
\qua
\sup_n  \|g_n(F_n)\|_{\DW}^2 <\infty,
\tio 
Lemma 1.2.3 \cite{nualart} 
states that this implies that $g(F) \in \DW$  and that  
$(D_{\cdot,0}\, g_n(F_n))$ converges to  $D_{\cdot,0}\,g(F)$ 
in the weak topology of $L_2(\Om; L_2(\R_+\times\{0\}))$. \\

Now  $ \E |g_n'(F_n)|^2 \le L^2_g $  implies that there exists a subsequence 
$(g_{n_k}'(F_{n_k}))_k$ which converges to some $G \in L_2$ in the 
weak topology  of $L_2.$ One can show that $|G| \le L_g$ a.s. Hence for 
any element $\al \in L_{\infty}(\Om; L_2(\R_+\times\{0\}))$ we have
\qua
&& \less \lim_{k \to \infty} \E \int_{\R_+}  
   g_{n_k}'(F_{n_k})(D_{t,0}\,F_{n_k}) \al(t) dt \\
&=&  \lim_{k \to \infty} \E \kla g_{n_k}'(F_{n_k}) \int_{\R_+}  
  (D_{\cdot,0}\,F_{n_k}) \al(t) dt \mer \\
&=& \lim_{k \to \infty} \E \kla g_{n_k}'(F_{n_k}) \int_{\R_+}  
  (D_{t,0}\, F) \al(t) dt \mer = \E \kla G  \int_{\R_+}  
  (D_{t,0}\, F) \al(t) dt \mer,
\tio
since $ | \E  g_{n_k}'(F_{n_k}) \int_{\R_+}  
  (D_{t,0} \,F_{n_k} - D_{t,0} \,F) \al(t)dt| \le L_g \| \al\|_{L_2(\Om; L_2(\R_+))}
        \|F_{n_k}-F\|_{\DW}  $
converges to zero for $k \to \infty$ and 
$\int_{\R_+} (D_{t,0} \,F) \al(t) dt \in L_2 $ because 
$\|\al\|_{L_2(\R_+\times\{0\})}$ is bounded. 
Consequently, 
\qua
\E \int_{\R_+}  D_{t,0}\, g(F)\al(t) dt =  \E  \int_{\R_+}G  
(D_{t,0} \,F) \al(t) dt 
\tio  
which implies $D_{t,0}\,g(F) =GD_{t,0} \, F$ \pl $ dt \otimes\mass -\text{a.e.}$  \\

(b) Let  $(F_n)_n \seq \mathcal{S}$ be a sequence such that
$\DN-\lim F_n =F.$ Since the expression
$ Z(t,x):= \frac{g(F + xD_{t,x}F)-g(F)}{x}\one_{\R_0}(x)$ is in 
$L_2(\m \otimes \mass)$  it is enough to show that  $(Dg_n(F_n)\one_{\R_0})$
converges in $ L_2(\m \otimes \mass)$ to $Z$ 
where $(g_n)$ is the sequence constructed in (a).  
Choose $T>0$ and $L>0$ large enough and $\delta >0$ sufficiently small such that
\qua
\limsup_n \E \int_{([0,T] \times \{\delta\le|x| \le L \})^c} |Z(t,x)|^2 
          + |D_{t,x} g_n(F_n)|^2 d\m(t,x) <\vare.
\tio
Then, for $n \ge n_0,$ 
\qua
&& \les  \| Z- Dg_n(F_n)\one_{\R_0}\|^2_{ L_2(\m \otimes \mass)} \\
&\le&  \vare \\
&&+ 2\E \int_{[0,T] \times \{\delta\le|x| \le L \}} |Z(t,x) -D_{t,x} g(F_n) |^2
  d\m(t,x) \\
&&+  8 \delta^{-2} T \mu(\{\delta \le |x| \le L \}) \|g-g_n\|^2_{\infty}.
\tio
Hence we obtain  \eqref{D_{t,x}-part}
from the uniform convergence of $g_n$ to $g$
and the Lipschitz continuity of $g.$ 
\boxx

\begin{proposition} Let $g:\R \to \R$ be Lipschitz continuous. 
Then $F \in \DD_{1,2}$ implies $g(F) \in \DD_{1,2},$ 
where $Dg(F)$ is given by \eqref{D_{t,0}-part} and 
\eqref{D_{t,x}-part}. 
\end{proposition}

\proof The assertion is an immediate consequence from Lemma  
\ref{Lipschitz-lemma} and \eqref{intersection}. \boxx



\end{document}